\documentclass{svproc}

\usepackage{amsmath}
\usepackage{mathtools}
\usepackage{amssymb}
\usepackage{tabularx}
\usepackage{xspace}
\usepackage{booktabs}
\usepackage{tikz}
\usepackage{paralist}
\usepackage{subcaption}
\usepackage{csquotes}

\usepackage{newtxtext}
\usepackage{newtxmath}
\usepackage{color}

\usetikzlibrary{arrows.meta,calc}

\usepackage[
  pdfencoding=auto,
  psdextra,
  colorlinks=true,
  linkcolor=blue,
  citecolor=blue,
  urlcolor=blue, 
  filecolor=blue,
  menucolor=blue,
  bookmarksopen=true,
  bookmarksnumbered=true,
  pdfstartview={FitH},
  pdfauthor={Ned Nedialkov and John Pryce},
]{hyperref}

\usepackage{doi}
\usepackage{cleveref}
\crefname{subsection}{subsection}{Subsections}
\Crefname{subsection}{Subsection}{Subsections}
\crefname{lstlisting}{Listing}{Listings}
\Crefname{lstlisting}{Listing}{Listings}
\crefname{appendix}{Appendix}{Appendices}
\Crefname{appendix}{Appendix}{Appendices}
\crefname{remark}{remark}{remarks}
\Crefname{remark}{Remark}{Remarks}

\newcommand\nc[1]{\newcommand{#1}}
\nc\rnc[1]{\renewcommand{#1}}

\nc\cmt[1]{{\color{red}#1}}

\nc\R{{\mathbb R}}
\nc\bmx[1]{\begin{bmatrix}#1\end{bmatrix}}
\nc\dom{\mathcal{M}}
\nc\pbg[1]{\bigl(#1\bigr)}
\nc\X{\!\times\!}
\nc\obsmap{\mathrm{\Phi}}
\nc\medbmx[1]{\mbox{\footnotesize$\bmx{#1}$}}
\nc\rfrac[2]{#1/#2} 
\nc\tc[2]{#1_{#2}}
\nc\tcf[2]{(#1)_{#2}}

\nc\tbeg{t_0}
\nc\tend{t_\text{end}}
\nc\xp{\dot{x}}
\nc\up{\dot{u}}
\nc\vp{\dot{v}}
\nc\tsexpand[1]{#1_0 + #1_1 s + #1_2 s^2 + \cdots}

\newcommand{\dbd}[2]{\frac{\partial #1}{\partial #2}}
\newcommand{\der}[2]{#1^{(#2)}}
\nc\adfk[1]{\operatorname{ad}_f^{#1}}
\nc\lider[1]{L_f^{#1}\,}
\rnc\d{\mathrm{d}}

\nc\TS{Taylor series\xspace}
\nc\TC{Taylor coefficient\xspace}
\nc\TCs{Taylor coefficients\xspace}

\nc\Lder{Lie derivative\xspace}
\nc\Lders{\Lder{s}\xspace}
\nc\LC{Lie coefficient\xspace}
\nc\LCs{\LC{}s\xspace}

\nc\matlab{MATLAB\xspace}
\nc\cl{code-list\xspace}
\nc\sode{sub-ODE\xspace}
\nc\sodes{sub-ODEs\xspace}
\nc\baos{BAOs\xspace}
\nc\Adtayl{\textsc{Adtayl}\xspace}
\nc\adtayl{\li{adtayl}\xspace}
\nc\odets{\li{odets}\xspace}

\usepackage{listings}
\usepackage{etoolbox}
\usepackage{xcolor}
\usepackage{stmaryrd}
\usepackage{eqparbox,calc}
\usepackage[nobreak=true]{mdframed}

\makeatletter
\newcommand{\currentfontsize}{\fontsize{\f@size}{\f@baselineskip}\selectfont}
\makeatother

\nc\li[1]{\text{\lstinline[basicstyle=\fontfamily{lmss}\fontseries{m}\selectfont]{#1}}}

\nc\para[1]{{\bf #1}}

\nc\HIr[1]{{\color{red}#1}}
\nc\HIb[1]{{\color{blue}#1}}
\nc\RT[2]{{\color{red}#1} {\color{blue}#2} }
\nc\term[1]{{\fontfamily{lmss}\fontseries{m}\slshape\selectfont#1}}

\makeatletter
\def\rf#1{(\@rf#1,.)}
\def\@rf#1,{\ref{eq:#1}\@ifnextchar . {\@endrf}{, \@rf}}
\def\@endrf.{}
\makeatother

\begin{document}

\title{High-Order Lie Derivatives from Taylor Series in the ADTAYL Package}

\author{
	Nedialko S.~Nedialkov\inst{1} \and
	John D.~Pryce\inst{2}
}
\authorrunning{Nedialko S.~Nedialkov \and John D.~Pryce}
\tocauthor{Nedialko S.~Nedialkov, John D.~Pryce}

\institute{
	McMaster University, Hamilton, Ontario, Canada\\
	\email{nedialk@mcmaster.ca}
	\and
	Cardiff University, Cardiff, United Kingdom\\
	\email{PryceJD1@cardiff.ac.uk}
}\maketitle

\setcounter{tocdepth}{2}  

\begin{abstract}
	High-order \Lders are essential in nonlinear systems analysis.
If done	symbolically, their evaluation becomes increasingly expensive as the order increases.
We present a compact and efficient numerical approach for computing
	\Lders of scalar, vector, and covector fields using the {\matlab} \Adtayl
	package.
The method exploits a fact noted by R\"obenack: that these derivatives coincide, up to factorial scaling, 
	 with the Taylor coefficients of expressions built from a Taylor expansion
about a trajectory point and, when required, the associated variational matrix.
Computational results for a gantry crane model demonstrate
	orders of magnitude speedups over symbolic evaluation using the
	\matlab Symbolic Math Toolbox.

\end{abstract}

\section{Introduction}
We have developed \Adtayl~\cite{adtaylUserGuide}, a numerical \matlab package
for Taylor series arithmetic.
This paper presents its application to the efficient evaluation of
high-order \Lders in nonlinear dynamical systems, a task for which symbolic
methods often become impractical at high derivative orders.

Such
 \Lders play an important role in nonlinear systems theory,
arising in output differentiation, observability and controllability
analysis, feedback linearization, and normal-form construction
\cite{isidori1995nonlinear,sontag1998}. Symbolic computation of Lie
derivatives is supported in the computer algebra system  Maple
\cite{mapleDifferentialGeometry}, and can also
be carried out using the \matlab Symbolic Math Toolbox and SymPy
\cite{meurer2017sympy}. However, symbolic approaches typically exhibit
expression swell, making high-order evaluation increasingly expensive as the
derivative order grows; see, e.g.,~\cite{robenack2004computation,Robenack2011LIEDRIVERSA}.

Numerical approaches based on automatic differentiation provide an effective
alternative.
Prior work by R\"obenack and co-authors showed that
high-order \Lders can be computed efficiently by exploiting the fact that
factorial-scaled Lie derivatives coincide with Taylor coefficients of
suitably constructed expressions involving the system flow and, when
required, the associated variational matrix \cite{robenack2005computation,robenack2024toward,robenack2004computation}. Their ADOL-C implementation achieved substantial speedups over symbolic
methods~\cite{Robenack2011LIEDRIVERSA}.

Using the mathematical identities introduced by R\"obenack \cite{robenack2005computation,robenack2004computation}, we produce implementations using Taylor series arithmetic, whose core computations reduce to two lines of MATLAB code each for Lie derivatives of scalar, vector, and covector fields respectively.
The same code applies unchanged to families of such fields.
This simplicity arises from the features and existing machinery of \Adtayl.

Numerical experiments on a gantry crane model illustrate the efficiency of the approach,  compared with evaluation using the \matlab Symbolic Math Toolbox.

The remainder of the paper is organized as follows.
Section~\ref{sc:liedef} defines Lie derivatives and states the computational problem.
Section~\ref{sc:overview} describes the \Adtayl functionality relevant to our implementation.
Section~\ref{sc:lder} presents the computation of high-order \Lders for scalar, vector, and covector fields, and compares our approach with the LIEDRIVERS~\cite{Robenack2011LIEDRIVERSA} implementation.
Performance results are given in Section~\ref{sc:perf}, and conclusions in Section~\ref{sc:concl}.


\section{Lie Derivatives: Definitions and Problem Statement}\label{sc:liedef}

\begin{definition}[\Lder]\label{df:lie_derivative}~\\
  Let $\dom \subseteq \R^n$ be an open set, and let $f:\dom\to\R^n$ be a smooth vector field
  defining the autonomous ordinary differential equation
  \begin{align}\label{eq:mainode}
    \xp = f(x).
  \end{align}

  \medskip
  \noindent\textbf{(a) Lie derivative of a scalar field.}
  Let $h:\dom\to\R$ be a smooth scalar field.
  The \Lder of $h$ along $f$ is the directional derivative $\lider{}h:\dom\to\R$
  defined by
  \begin{align}\label{eq:Lmap}
    \lider{} h(x) = h'(x)\,f(x),
  \end{align}
  where $h'(x)\in\R^{1\times n}$ is the row gradient of $h$.

  \smallskip
  \noindent\textbf{(b) Lie derivative of a vector field.}
  Let $g:\dom\to\R^n$ be a smooth vector field.
  The \Lder of $g$ along $f$ is the \term{Lie bracket}
  \begin{align}\label{eq:Lbracket}
    \lider{} g(x) = [f,g](x) = g'(x)\,f(x) - f'(x)\,g(x),
  \end{align}
  where $f'(x),g'(x)\in\R^{n\times n}$ are Jacobians. 
  
  \smallskip
  \noindent\textbf{(c) Lie derivative of a covector field.}
  Let $\omega:\dom\to(\R^n)^*$ be a smooth covector field.  The Lie derivative of $\omega$ along $f$ is defined by
  \begin{align}\label{eq:Lcov}
    \lider{} \omega(x)
    = \bigl[\omega'(x)f(x)\bigr]^\top + \omega(x)f'(x), \quad\text{where}\quad
    \omega'(x) = \dbd{}{x}\,\omega(x)^\top \in \R^{n\times n}.
  \end{align}
\end{definition}
\smallskip
  
For a given $\dom$ and vector field $f$ on $\dom$, the mapping $X \mapsto \lider{}X$ is an operator of the space of smooth fields of the given type (scalar $h$, vector $g$, or covector $\omega$) into itself.
It is linear: for all such $X,Y$ and scalars $\alpha,\beta$,
\[
  \lider{}(\alpha X+\beta Y)=\alpha\,\lider{}X+\beta\,\lider{}Y.
\]
The \term{iterated Lie derivative} $\lider{k}$ is this operator's $k$th power: with $X$ denoting $h$, $g$, or $\omega$,
  \begin{align}\label{eq:Lderit}
    \lider{0}X = X, \qquad
    \lider{k} X = \lider{} \lider{k-1} X,
    \quad \text{for } k \ge 1.
  \end{align}
It is common to denote $\lider{k}$ of a vector field by $\adfk{k}$, the \term{adjoint operator}.  
\medskip

\noindent\textbf{Problem statement.}
Given a smooth vector field
$f:\dom\to\R^n$, a smooth field $X$ on $\dom$ (scalar, vector, or covector),
an initial state $x_0\in\dom$, and an integer $p\ge 1$, compute the iterated
Lie derivatives
\begin{align}\label{eq:problem}
  &\lider{k} X(x_0), \quad k = 0, 1, 2, \dots, p.
\shortintertext{We call}
  \frac{1}{k!}\,&\lider{k}X(x_0), \qquad k \ge 0, \notag
\end{align}
the \term{Lie coefficients} of $X$ at $x_0$ along $f$.
Since these coefficients arise naturally in our theory, we compute Lie
coefficients; the actual derivatives $L_f^k X(x_0) $ are then obtained by multiplying by $k! $.

\section{\Adtayl Overview}
\label{sc:overview}

\Adtayl, \underline{A}utomatic \underline{D}if\-ferentiation by
\underline{Tayl}or Series, is a purely numerical \matlab package with two main
components:
\begin{compactitem}
	\item \adtayl, a class for univariate Taylor series arithmetic, and
	\item \odets, a solver for Taylor-based integration of ordinary
	differential equation (ODE) initial-value problems (IVPs).
\end{compactitem}
Together, \adtayl and \odets provide the machinery used here to evaluate \LCs.

\Cref{ss:adtayl} describes the features of \adtayl that are relevant to the
implementation in \Cref{sc:lder}.
\Cref{ss:odets} presents the corresponding features of \odets.
\Cref{ss:tcs} summarizes the method used to compute \TCs in both components.
A complete description of this package is provided in the \Adtayl User
Guide~\cite{adtaylUserGuide}.

\subsection{The \adtayl Class}\label{ss:adtayl}

The \adtayl class implements \matlab-style arrays whose elements are
truncated Taylor series  in a single independent variable.
Concretely, each entry in an \adtayl array stores the truncated Taylor
expansion of a scalar function (denoted here by $u(t)$) about some $t=t_0$ to a
user-specified order $p\ge0$,
\begin{equation*}
	u(t_0+s)=u_0+u_1 s+\cdots+u_p s^p+O\bigl(s^{p+1}\bigr),
	\qquad
	u_k=\frac{\der{u}{k}(t_0)}{k!}.
\end{equation*}
That is, the entry corresponding to $u(t)$ stores the coefficients
$u_0, u_1, \ldots, u_p$. All entries in an \adtayl array store series of the same order.

This class allows \matlab expressions to be evaluated without modifying user code when some or all operands are \adtayl objects.
Expressions may freely mix \adtayl and numeric arrays of compatible sizes, provided that all \adtayl operands have the same order, or that one operand has order zero; numeric arrays are interpreted as Taylor series of order zero.
Binary operations involving \adtayl objects follow the standard \matlab rules for element-wise operations on compatible array sizes.

This functionality is implemented by overloading elementwise arithmetic operators
(\li{plus}, \li{minus}, \li{times},...) and a full set of elementary functions
(e.g., \li{exp}, \li{log}, \li{sin}, \li{asin}, \ldots), which all operate on
truncated Taylor series and discard terms of order $>p$.
In addition, \adtayl provides access and array-manipulation operations (reshaping,
concatenation, indexing, slicing) that act on the stored coefficients.

The following two matrix-level
operations are
used in \Cref{sc:lder}.

Let
\li{A} and \li{B} denote \adtayl objects representing matrix-valued
functions $A(t)$ and $B(t)$. \begin{compactitem}\setlength{\itemsep}{3pt}
	\item The matrix product \li{A*B} computes the truncated Taylor expansion
	of the product $A(t)B(t)$.
	\item The backslash operator \li{A\\B} computes the truncated Taylor
	expansion of $X(t)$ such that
	$A(t)\,X(t)=B(t)$, provided that the constant term of $A(t)$ is a
	nonsingular matrix.
\end{compactitem}
In both cases, \li{A} and \li{B}  must have compatible array dimensions and orders.

By default, \adtayl operates in double precision. It also works with  the
variable-precision arithmetic (\li{vpa}) of \matlab. In both double and
\li{vpa} modes, \adtayl supports real and complex evaluation of all
standard functions.

\subsection{The \odets Solver}\label{ss:odets}
The \odets solver is designed to integrate initial-value problems of the
form
\begin{align}\label{eq:mainivp}
	\dot{x} = F(t,x), \quad x(t_0) = x_0, 
\end{align}
where $F$ is
built from the
four basic arithmetic operations and suitable standard functions. We exclude
here non-differentiable functions like \li{abs}, \li{min} and \li{max},  and conditional
structures such as \enquote{if} statements.

Given a \matlab function that evaluates the right-hand side $F$ of the
ODE, \odets preprocesses the function  using operator overloading to construct
a \cl object in memory, which encodes the evaluation of $F$
\cite{nedialkovpryce2025}. This preprocessing step is
conceptually similar to the taping mechanism in ADOL-C~\cite{adolc}.

During each integration step, \odets interprets this \cl to compute
the Taylor coefficients of the solution at the current time point. It
then advances the solution by evaluating the resulting truncated
Taylor series with a suitable step size.
\smallskip

The \odets solver also accepts IVP ODEs in autonomous form
\begin{align}\label{eq:ivp}
	\dot{x} = f(x), \qquad x(t_0) = x_0,
\end{align}
which is the case of interest here. In addition to computing the Taylor coefficients of the solution
$x(t)$ to \rf{ivp}, \odets can also compute the Taylor coefficients of
$J(t) = \partial x(t)/\partial x_0$, which satisfies the variational equation
\begin{equation}\label{eq:var}
	\dot{J}(t) = f'\bigl(x(t)\bigr)\, J(t),
	\qquad J(t_0) = I,
\end{equation}
where $I$ is the identity matrix.

\subsubsection{The \li{taylcoeffs} function.}
For the above computation, we provide the function
\begin{lstlisting}[numbers=none]
[x,J] = taylcoeffs(f,x0,p)
\end{lstlisting}
Argument \li{f} is either a function handle for evaluating  $f$ or a \cl pre-constructed from $f$; argument \li{x0} holds initial state vector $x_0$; and \li{p} is the desired expansion order~$p$.

Internally, \li{taylcoeffs} uses the \cl infrastructure of \odets to compute
the Taylor coefficients of the solution $x(t)$ to \rf{ivp}
 to order $p$, and stores them in the \adtayl object \li{x}. When the output
\li{J} is requested, it computes the Taylor coefficients of  $J(t)$ simultaneously and returns them in the \adtayl object \li{J}.

\subsection{Computing Taylor Coefficients}\label{ss:tcs}

The computation of Taylor coefficients in \Adtayl is based on five recurrences: four for
the basic arithmetic operations (\baos), $+,-,\times,/$, and one for the
\sode\ operation~$\odot$, described below. We summarize these here. For details on the \sode method, see~\cite{nedialkovpryce2025}.

\subsubsection{Recurrences for the \baos.}

Let the Taylor coefficients of $u(t)$ and $v(t)$ be known up
to order $k$. For the \baos, the $k$th Taylor coefficients are given by
\begin{align}
	\tcf{u \pm v}{k} & = \tc{u}{k} \pm \tc{v}{k},
	\label{eq:tca}                                              \\[2pt]
	\tcf{u v}{k}     & = \sum_{i=0}^{k} \tc{u}{i}\,\tc{v}{k-i},
	\label{eq:tcm}                                              \\[2pt]
	\tcf{u/v}{k}     & = \frac{1}{\tc{v}{0}}
	\Bigl(\tc{u}{k}
	- \sum_{i=0}^{k-1} \tc{v}{k-i}\,\tcf{u/v}{i}\Bigr),
	\qquad \tc{v}{0} \neq 0.
	\label{eq:tcd}
\end{align}
In \rf{tcd}, the coefficient $\tcf{u/v}{k}$ depends on $\tc{u}{k}$ and on coefficients of strictly lower order.

Assuming the solution of \rf{ivp} admits a Taylor expansion about $t=t_0$, the $(k+1)$st Taylor coefficient of $x(t)$ satisfies
\begin{equation}\label{eq:xkp1}
  x_{k+1} = \frac{1}{k+1}\,f_k,
\end{equation}
where $x_{k+1}=\tfrac{\der{x}{k+1}(t_0)}{(k+1)!}$ and $f_k=\frac{1}{k!}\frac{\d^k}{\d t^k}f\bigl(x(t)\bigr)\big|_{t=t_0}$\;.

\smallskip
If $f$ is composed solely of \baos, then \rf{xkp1} together with
\rf{tca}, \rf{tcm}, and \rf{tcd} yields a complete procedure for computing the coefficients
$x_1, x_2, \ldots$.

\subsubsection{Standard functions via \sodes.}

We begin with an example and then describe the general scheme.

\begin{example}
	Suppose the Taylor coefficients of $u(t)$ are known up to order $k$, and we wish
	to compute the coefficients of $v(t)=\exp\big(u(t)\big)$ to the same order. The
	exponential satisfies
	\[
		\frac{\d v}{\d u} = v,
	\]
	and thus satisfies the ODE
	\begin{align}\label{eq:expsubode}
		\vp = v\up.
	\end{align}

	Applying \rf{xkp1} and  \rf{tcm} to  \rf{expsubode}, the $k$th coefficient for $k\ge 1$ is 
	\begin{align*}
		v_k
= \frac{1}{k}(v\dot u)_{k-1}
= \frac{1}{k}\sum_{i=1}^{k} i\,u_i\,v_{k-i}.	\end{align*}
	This recurrence is started with  $v_0 = \exp(u_0)$, which can be viewed as the initial condition at $t=0$ for \rf{expsubode}.
\end{example}

\nc\subo{\varphi}
\nc\stdf{\psi}

In general, let $\stdf(u)$ be an $m$-vector function of a scalar argument and
set $v=\stdf(u)$. If the derivative $\stdf'(u)\in\R^m$, with respect to the scalar $u$, can be written in the
form
\begin{align}\label{eq:huv0}
  \stdf'(u) = \subo\bigl(u,\stdf(u)\bigr),
\end{align}
then $v(t)$ satisfies the sub-ODE
\begin{align}\label{eq:gsubode}
  \dot v = \subo(u,v)\,\dot u.
\end{align}

\begin{example}
\nc\cs{\text{cs}}
\nc\smallbmx[1]{\left[\begin{smallmatrix}#1\end{smallmatrix}\right]}
	\rnc{\_}{}
	Examples of \sodes for standard functions:
	
	{
	\noindent
	\makebox[0pt][l]{
		\hspace{-1cm}
		\begin{minipage}{\linewidth}
			\begin{align*}
				\begin{tabular}{rl@{\qquad has }l@{\qquad $\subo(u,v)=$ }l}
					(i)   & $v = \exp(u)$             & $\vp = \_v \up$,        & $v$    \\
					[1ex]
					(ii)  & $v = u^c$, ($c$ constant) & $\vp = (c\_v/\_u) \up$, & $cv/u$ \\
					[2ex]
					\multicolumn{4}{l}{\text{and, defining $\cs(u) = \smallbmx{\cos(u)    \\ \sin(u)}$,}} \\
					[2ex]
					(iii) & $v  = \bmx{v_1                                                 \\ v_2} = \cs(u) $ &
					$\vp = \bmx{-\_v_2                                                    \\ \_v_1} \up$, &
					$\bmx{-v_2                                                            \\ v_1}$
				\end{tabular}
			\end{align*}
		\end{minipage}
	}}

\end{example}

\smallskip
In practice,  $\subo(u,v)$ contains  only  \baos. 
 For a
standard function $\stdf(u)$ with \sode \rf{gsubode}, the coefficients are computed using the $\odot$ operation
\begin{align}\label{eq:ddot}
	v_{k}  =  (\subo\odot_{\stdf} u)_{k} =
	\begin{cases}
		\stdf(u_0),                                                 & k=0, \\[1ex]
		\displaystyle\frac{1}{k}\sum_{i=1}^{k} i\,u_i\,\subo_{k-i}, & k>0.
	\end{cases}
\end{align}
Here $\subo_j$ denotes the $j$th Taylor coefficient of $\subo\bigl(u(t),v(t)\bigr)$.

\begin{example}
For the case (ii) in the previous example, using \rf{ddot} the coefficients are computed as 
\begin{align*}
	v_{k}  =  
	\begin{cases}
		{u_0}^c,                                                 & k=0, \\
		\frac{1}{k}\sum_{i=1}^{k} i\,u_i\,(cv/u)_{k-i}, & k>0,
	\end{cases}
\end{align*}
where the  $(cv/u)_{k-i}$ are found using \rf{tcd}.
\end{example}

\begin{remark}
With the \sode machinery in place, implementing standard functions is remarkably concise.
In about 130 lines of \matlab code in \odets---about 170 lines in \adtayl---we have implemented the exponential, logarithm and power functions, 
all direct and inverse trigonometric and hyperbolic functions;
root functions \li{sqrt}, \li{nthroot}, and functions such as \li{expm1}, \li{log1p}, and \li{atan2}.
\end{remark}

\begin{remark}
	Our approach differs from that of established
	systems such as
	TIDES~\cite{abad2012algorithm},
	ATOMFT~\cite{Chang1994a},
	ADOL-C~\cite{adolc},
	TAYLOR~\cite{jorba2005software},
	and
	FADBAD++~\cite{FADBAD++}.
	In those packages, separate coefficient recurrences are  implemented
	for each standard function.
	
\end{remark}

\subsubsection{Jacobians of \TCs.}
Write 
\[
  x(t) = \sum_{k=0}^{\infty} x_k\, t^k .
\]
Each $x_k$ is a function of the initial $x_0$. 
Differentiating with respect to this yields
\begin{align}\label{eq:Jk}
  J(t) = \frac{\partial x(t)}{\partial x_0}
  = \sum_{k=0}^{\infty} \frac{\partial x_k}{\partial x_0}\, t^k .
\end{align}
To compute the coefficients in \rf{Jk}, we augment 
recurrences \rf{tca}, \rf{tcm}, \rf{tcd}, and \rf{ddot} to propagate 
gradients with respect to $x_0$. For example,  \rf{tcm} is augmented by
\begin{align*}
  \nabla\tcf{u v}{k}
  = \sum_{i=0}^{k}
  \Bigl(
   \nabla\tc{u}{i}\,\tc{v}{k-i}+ \tc{u}{i}\,\nabla\tc{v}{k-i}
  \Bigr).
\end{align*}
The input $i$th gradient is initialized as the $i$th unit vector. 

\section{Computing Lie Coefficients}\label{sc:lder}

The following three subsections show how we compute the \LCs in \rf{problem}
for scalar, vector, and covector fields. \Cref{ss:families} shows that, in
each case, the \Adtayl implementation extends without modification to families
of such fields. \Cref{ss:robenack} highlights the differences from the
approach in~\cite{Robenack2011LIEDRIVERSA}.

In the code snippets below, the resulting \LCs are returned in an \adtayl object \li{L}. Their
numerical values are extracted as \li{L.gettc(k)}, for $k = 0,1,\ldots,p$.

\subsection{Scalar Field}\label{ss:map}
Here $h$ is a scalar field.
Let  $y(t) = h\big(x(t)\big)$, where $x(t)$ satisfies \rf{ivp}. Then 
\begin{align}\label{eq:ykt}
	y^{(k)}(t) = \lider{k}h\big(x(t)\big),
\end{align}
and therefore the $k$th \TC of $y(t)$ at $t=0$ is
\begin{align}\label{eq:yklieder}
	y_k = \frac{y^{(k)}(0)}{k!} = \frac{1}{k!}\frac{d^k}{dt^k} h\big(x(t)\big)\Big|_{t = 0}
	= \frac{1}{k!} \lider{k}h(x_0).
\end{align}
Hence the $k$th \LC at $x_0$ is precisely the $k$th \TC of  $h\pbg{x(t)}$ at $t=0$; see~\cite{robenack2024toward}.
The following \Adtayl code computes the $y_k$ in \rf{yklieder}:
\begin{lstlisting}[numbers=none,caption={\LCs of a scalar field.},label=ls:liesc]
x = taylcoeffs(f,x0,p);
L = h(x);
\end{lstlisting}

\begin{compactitem}
	\item Here, \li{f} is a function handle or a \cl object for evaluating
	the vector field $f$.

	\item \li{h} is a function handle for evaluating $h$.

	\item The \adtayl object \li{x} contains $x_0, x_1, \ldots, x_p$.

	\item Thus \li{h(x)} evaluates $h(x)$ in \adtayl arithmetic and propagates \li{x}
	through the expressions defining $h$. The resulting \adtayl object  \li{L} contains the coefficients $y_0, \ldots, y_p$ of
	$y(t) = h\bigl(x(t)\bigr)$, which by \rf{yklieder} are the desired \LCs.

\end{compactitem}

\subsection{Vector Field}\label{ss:lievec}
Here $g$ is a vector field. 
Let $Z(t)$ be defined by
\begin{align}\label{eq:Z}
	\dot{Z}(t) = -Z(t)\, f'\pbg{x(t)}, \qquad Z(0) = I.
\end{align}
It is easily seen that $Z(t)=J^{-1}(t)$ as a formal Taylor series (since
$J(0)=I$), where $J(t)$ satisfies the variational equation \rf{var}.
Using the derivation in ~\cite[Eq.~(25)]{robenack2004computation}, one has
\begin{align}
	\sum_{k=0}^{\infty} \frac{1}{k!}\, \adfk{k} g(x_0)\,t^k & = Z(t)\, g\pbg{x(t)} \notag              \\
	                                                        & = J^{-1}(t)\, g\pbg{x(t)}. \label{eq:Zt}
\end{align}
Hence,
\[
	\text{\LC $\frac{1}{k!} \adfk{k} g(x_0)$ is the $k$th \TC  of $J^{-1}(t)\, g\pbg{x(t)}$ at $t=0$} .
\]

Computing these in \Adtayl is done by
\begin{lstlisting}[numbers=none,caption={Lie coefficients of  a vector field.},label={ls:lievec}]
[x,J] = taylcoeffs(f,x0,p);
L = J\g(x);
\end{lstlisting}

\begin{compactitem}
	\setlength{\itemsep}{-2pt}

	\item The \adtayl objects \li{x} and \li{J} contain the  \TCs  of $x(t)$ and $J(t)$ to order $p$, respectively.

	\item \li{g(x)} is a function handle that returns an \adtayl object containing the \TCs of $g\big(x(t)\big)$ to order $p$.

	\item Since $J(0)=I$ is nonsingular, power series $J(t)$ is invertible, and \li{L=J\\g(x)} computes the \TCs of $J^{-1}(t)\, g\pbg{x(t)}$ to order $p$, returned in \li{L}.
	
\end{compactitem}

\subsection{Covector Field}
Here $\omega(x)$ is a covector field.
Using the result in \cite[Sec.~5.2]{robenack2005computation}, one has
\begin{align}\label{eq:liecov}
	\sum_{k=0}^{\infty} \frac{1}{k!}\, \lider{k} \omega(x_0)\, t^k
	= \omega\pbg{x(t)}\, J(t),
\end{align}
where $J(t)$ satisfies \rf{var}.
Hence,
\[
	\text{\LC $\frac{1}{k!}\,\lider{k}\omega(x_0)$ is the $k$th \TC of $\omega\pbg{x(t)}\,J(t)$} .
\]

Computing these in \Adtayl is done by
\begin{lstlisting}[numbers=none,
                   captionpos=b, label=ls:liecov,
                   caption={\protect\LCs of a covector field},
                   abovecaptionskip=10pt]
[x,J] = taylcoeffs(f,x0,p);
L = w(x)*J;
\end{lstlisting}
Here \li{w(x)}   returns a row vector (covector), so that \li{L} represents $\omega\big(x(t)\big )\,J(t)$.

\subsection{Families of Scalar, Vector, and Covector Fields.}\label{ss:families}

Without modification, the Listings \ref{ls:liesc}, \ref{ls:lievec}, \ref{ls:liecov} handle families of fields.
Since a vector field is implemented as a column vector, \adtayl supports a 1D array of them stacked side by side into a matrix.
Similarly, since a covector field is a row vector, \adtayl supports a 1D array of them stacked vertically into a matrix.
Scalar fields can be in an arbitrarily-shaped array.

\subsubsection{Scalar fields.}
An \adtayl array can have any number of dimensions. In general, \li{h(x)} in Listing~\ref{ls:liesc}
can return an array of size
$[s_1,  s_2,\, \ldots,  s_l]$,
representing the result of evaluating an array $h(x)$ of $s_1 \X s_2 \X \cdots \X s_l$ separate scalar fields.
Listing~\ref{ls:liesc} remains valid and computes the Lie coefficients for \emph{all} these scalar fields simultaneously, and
\[
	\frac1{k!}\lider{k} h(x_0) \quad\text{is the\  $s_1 \X s_2 \X \cdots \X s_l$ \ array \li{L.gettc(k)}}.
\]

Of interest in control theory is when this $h$ is an output map, generally regarded as a column vector:
\begin{definition}[Output map]
	Let $\dom$ be a smooth manifold representing the state space of a dynamical system.
	An \term{output map} is a smooth function
	\[
		h : \dom \to \R^m
	\]
	that assigns to each state $x \in \dom$ a vector $y = h(x) \in \R^m$ of observable quantities.
	Vector $y$ is called the system's \term{output}; its components $h_1, \dots, h_m$ are scalar fields on~$\dom$.
\end{definition}
So in this case we have componentwise
\begin{align*}
	\lider{k}h(x) = \bigl[\,\lider{k}h_1(x),\ \lider{k}h_2(x),\ \ldots,\ \lider{k} h_m(x)\,\bigr]^\top \in \R^{ m}.
\end{align*}

\subsubsection{Vector fields.}
Suppose we have $m$ smooth vector fields $g_1(x),\ g_2(x),\ \ldots,\ g_m(x)$.
Let
\[
	g(x) = \bigl[\,g_1(x),\ g_2(x),\ \ldots,\ g_m(x)\,\bigr] \in \R^{n\times m}.
\]
We define $\adfk{k}$ on this matrix column-wise:
\[
	\adfk{k} g(x)
	=
	\bigl[
		\adfk{k} g_1(x),\
		\adfk{k} g_2(x),\
		\ldots,\
		\adfk{k} g_m(x)
		\bigr],
	\qquad k \ge 0.
\]
Equation~\rf{Zt} remains true for such $g$.
The computation in \cref{ls:lievec} stays valid, and
\[
	\frac{1}{k!}\,\adfk{k} g(x_0) \quad\text{is the \  $n\X m$ matrix}\  \li{L.gettc(k)}.
\]

\subsubsection{Covector fields.}
Assume we have $m$ smooth covector fields $\omega_i$, $i=1, \ldots, m$. Let
\[
	\omega(x)
	=
	\begin{bmatrix}
		\omega_1(x) \\[4pt]
		\omega_2(x) \\[4pt]
		\vdots      \\[4pt]
		\omega_m(x)
	\end{bmatrix}
	\in \mathbb{R}^{m\times n},
\]
where each $\omega_i(x)$ is a row vector.
We define $\lider{k}$ on this matrix row-wise:
\[
	\lider{k}\omega(x)
	=
	\begin{bmatrix}
		\lider{k}\omega_1(x) \\[4pt]
		\lider{k}\omega_2(x) \\[4pt]
		\vdots               \\[4pt]
		\lider{k}\omega_m(x)
	\end{bmatrix},
	\qquad k \ge 0.
\]
Equation \rf{liecov} remains true for such $\omega$.
The computation in \cref{ls:liecov} stays valid, and
\[
	\frac1{k!}\lider{k} \omega(x_0) \quad\text{is the \  $m\X n$ matrix \ \li{L.gettc(k)}}.
\]

The remarkably concise \Adtayl computations described above are summarized 
\begin{table}[ht]
	\centering
		\begin{tabular}{l@{\hspace{1.5em}}l@{\hspace{2em}}l}
		\toprule
		Field type & Size     & Expression  \\ \midrule
		scalar     & any      & \li{h(x)}   \\
		vector     & $n \X m$ & \li{J\\g(x)} \\
		covector   & $m \X n$ & \li{w(x)*J} \\
		\bottomrule
	\end{tabular}
\end{table}

\subsection{Comparison with  LIEDRIVERS}\label{ss:robenack}

We summarize the LIEDRIVERS~\cite{Robenack2011LIEDRIVERSA} implementation using ADOL-C and note two  key differences from \Adtayl; see also \cite{robenack2005computation,robenack2024toward,robenack2004computation}.

LIEDRIVERS obtains the coefficients of $x(t)$ from the ADOL-C tape for $f$ by a forward pass through the tape.
The field $X$ is also taped, and the coefficients of $X\bigl(x(t)\bigr)$ are computed by a forward pass through
this tape, with all operations carried out on the truncated Taylor series for $x(t)$.

\begin{compactenum}[(1)]
\item A key difference in our approach is that we do not construct a \cl object for $X$; that is, $X$ is \enquote{not taped}.
Instead, $X\bigl(x(t)\bigr)$ is evaluated directly in \adtayl arithmetic. 

\item A second difference lies in the way Jacobian coefficients are computed.
In the vector and covector cases, LIEDRIVERS uses a subsequent reverse pass through the tape for $f$ to compute the coefficients $A_i$ of the Jacobian matrix $f'\bigl(x(t)\bigr)$,
\begin{align*}
	f'\bigl(x(t)\bigr)
	=
	A_0 + A_1 t + \cdots + A_k t^k + \cdots ,
\end{align*}
by reverse-mode automatic differentiation.
These coefficients are then used in the recurrences \rf{Zk} and \rf{Jk_adolc}.
By contrast we use only forward mode AD.
\end{compactenum}

\subsubsection{Vector field.}

LIEDRIVERS computes the coefficients of $Z(t)\,g\pbg{x(t)}$ in \rf{Zt}. The
coefficients of $Z(t)$ are formed using
\begin{align}\label{eq:Zk}
	Z_{k+1}
	=
	-\frac{1}{k+1}
	\sum_{i=0}^{k} Z_i A_{k-i},
	\qquad Z_0 = I.
\end{align}
This recurrence follows directly from \rf{xkp1} and \rf{tcm} when applied to
\rf{Z}. Writing
\[
	g\bigl(x(t)\bigr) = g_0 + g_1 t + \cdots + g_k t^k + \cdots ,
\]
the $k$th Taylor coefficient of $Z(t)\,g\bigl(x(t)\bigr)$ is obtained as
$\sum_{i=0}^{k} Z_i\, g_{k-i}$, which yields
$\frac{1}{k!}\,\adfk{k} g(x_0)$.

\subsubsection{Covector field.}

LIEDRIVERS computes the coefficients of $J(t)$ using
\begin{align}\label{eq:Jk_adolc}
	J_{k+1}
	=
	\frac{1}{k+1}\sum_{i=0}^{k} A_{k-i}\,J_i,
	\qquad J_0 = I.
\end{align}
This also follows directly from \rf{xkp1} and \rf{tcm} when
applied to \rf{var}. Similarly, writing
\[
	\omega\bigl(x(t)\bigr)
	=
	\omega_0 + \omega_1 t + \cdots + \omega_k t^k + \cdots ,
\]
one finds the $k$th coefficient of $\omega\pbg{x(t)}\,J(t)$ as $\sum_{i=0}^{k} \omega_i J_{k-i}$.

\section{Performance}\label{sc:perf}
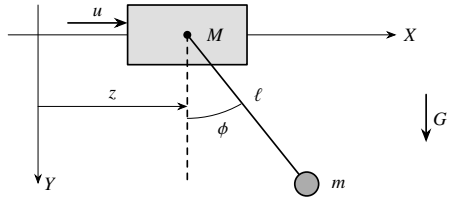
\begin{figure}[ht]
	\centering
	\resizebox{6.cm}{!}{

\begin{tikzpicture}[>=Stealth, line cap=round, line join=round, thick]

\coordinate (O) at (0,0);
\coordinate (Ox) at (6, 0);
\coordinate (Oy) at (0, -2.5);
\coordinate (CoM) at (2.5,0);
\coordinate (pend) at (4.5,-2.5);

\draw[->, thin] (O)+(-.5,0) -- (Ox) node[right] {$X$};
\draw[->, thin] (O) +(0,.5)-- (Oy) node[right] {$Y$};

\coordinate (g) at ($(Ox)+(.5,-1)$);
\draw[->] (g) -- ($(g)+(0,-.8)$) node[midway, right] {$G$};

\def\cartW{2}  
\def\cartH{1}  

\coordinate (CartSW) at ($(CoM)+(-0.5*\cartW,  0.5*\cartH)$); 
\coordinate (CartNE) at ($(CoM)+( 0.5*\cartW, -0.5*\cartH)$); 
\draw[fill=gray!25] (CartSW) rectangle (CartNE);
\node[right] at ($(CoM)+(.2,0)$) {$M$};

\coordinate (u) at ($(O)+(.5,0.2)$);
\draw[->] (u) -- ($(CoM)+(-0.5*\cartW,  .2)$) node[midway, above] {$u$};

\filldraw[black] (CoM) circle (1.5pt);
\draw (CoM) -- (pend);
\filldraw[fill=gray!65] (pend) circle (0.20);
\node[right] at ($(pend)+(0.30,0)$) {$m$};
\node[right] at ($(pend)+(-1,1.5)$) {$\ell$};

\draw[dashed] (CoM) -- ($(CoM |- pend)$) ;

\draw[->,thin] ($(O)+(0,-1.2)$) -- ($(CoM) +(0,-1.2)$) node[midway, above] {$z$};

\draw[thin] ($(CoM)+(0,-1.4)$)
  arc[start angle=-90, end angle=-60, radius=1.8];
\node at ($(CoM)+(0.6,-1.6)$) {$\phi$};
\end{tikzpicture}}
	\caption{Schematic of the gantry crane model.\label{fig:schem}}
\end{figure}

We report the performance of our approach 
on  the gantry crane model from
\cite{Robenack2011LIEDRIVERSA}, \Cref{fig:schem}. The system consists of
a cart of mass $M$ moving horizontally and carrying a load of mass $m$,
suspended by a massless cable of length $\ell$ at an angle $\phi$. The
input $u(t)$ is a horizontal force applied to the cart, and gravity $G$
acts vertically downward.

The equations of motion, where $z$ is horizontal position, are
\begin{align}\label{eq:motion}
	\begin{pmatrix}
		m + M           & m \ell \cos\phi \\
		m \ell \cos\phi & m \ell^2
	\end{pmatrix}
	\begin{pmatrix}
		\ddot z \\
		\ddot \phi
	\end{pmatrix}
	+
	\begin{pmatrix}
		- m \ell \dot\phi^2 \sin\phi \\
		 m \ell G  \sin\phi
	\end{pmatrix}
	=
	\begin{pmatrix}
		u \\
		0
	\end{pmatrix}.
\end{align}
With the state vector $x=(z,\phi,\dot z,\dot\phi)$, \rf{motion}
as a first-order system is of the form
\begin{align}\label{eq:gantry}
  \dot x = f(x) + g(x)u.
\end{align}
The expressions for the vector fields $f$ and $g$, and for the output map  $y = h(x)$ for the load position are given in Appendix~\ref{sc:gantryeqns}.
From the control-theory viewpoint, this has the frequently occurring {\em control-affine} form, meaning the input $u$ appears only linearly.
\smallskip

The tests consisted in computing the iterated Lie brackets $\adfk{k} g(x_0)$ and \Lders $\lider{k} h(x_0)$ for orders $k = 1,\ldots,10$.
We used the parameter values $M=1$, $m=1$, $\ell=1$, and $G=9.81$, and the initial state $x_0 = (1,0.2, -0.5,-0.4)^\top$, as given in the  \li{GantryCrane.cpp} file from the ADOL-C distribution.

We computed the derivatives (a) with the \Adtayl package, (b) with the MATLAB Symbolic Math Toolbox, evaluating the resulting symbolic expressions, (c) by re-running \li{GantryCrane.cpp} with ADOL-C.
The three sets of numbers were very close, giving us confidence in \Adtayl's algorithms.

\smallskip

\nc\scl{.48}

The tests were done on a Mac Studio (Apple M4 Max, 64\,GB RAM) with \matlab R2025b.
The timings reported below are wall-clock times.
For the symbolic approach, they include both the construction of the symbolic expressions  and their subsequent numerical evaluation.
For \Adtayl, they include the generation of the \cl for the vector field $f$ and the computation of the corresponding derivatives.

\begin{figure}[htb]
  \centering
  \begin{subfigure}{\scl\linewidth}
    \centering
    \includegraphics[width=\linewidth]{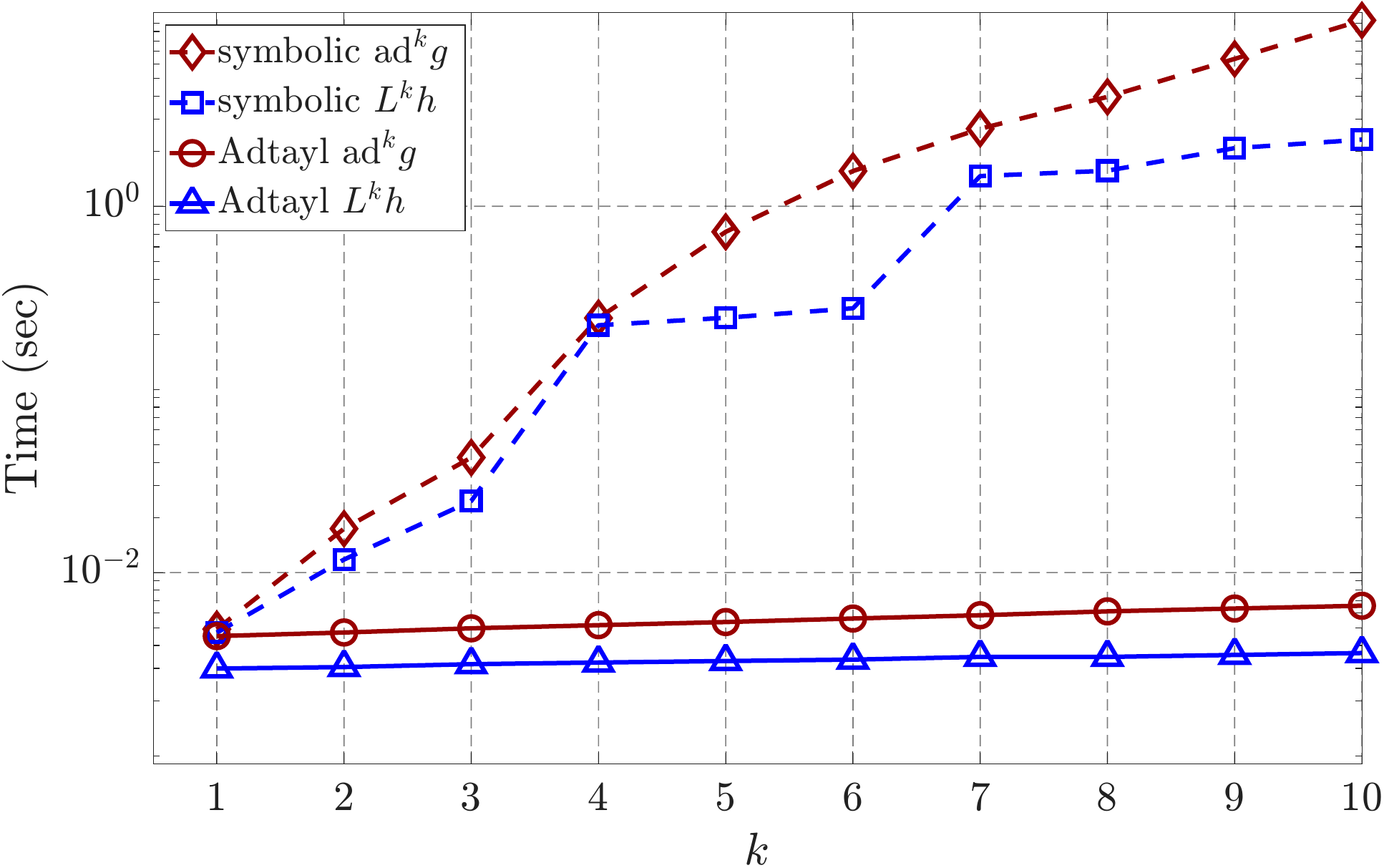}
    \caption{Time versus order.}
    \label{fig:combined}
  \end{subfigure}\hfill
  \begin{subfigure}{\scl\linewidth}
    \centering
    \includegraphics[width=\linewidth]{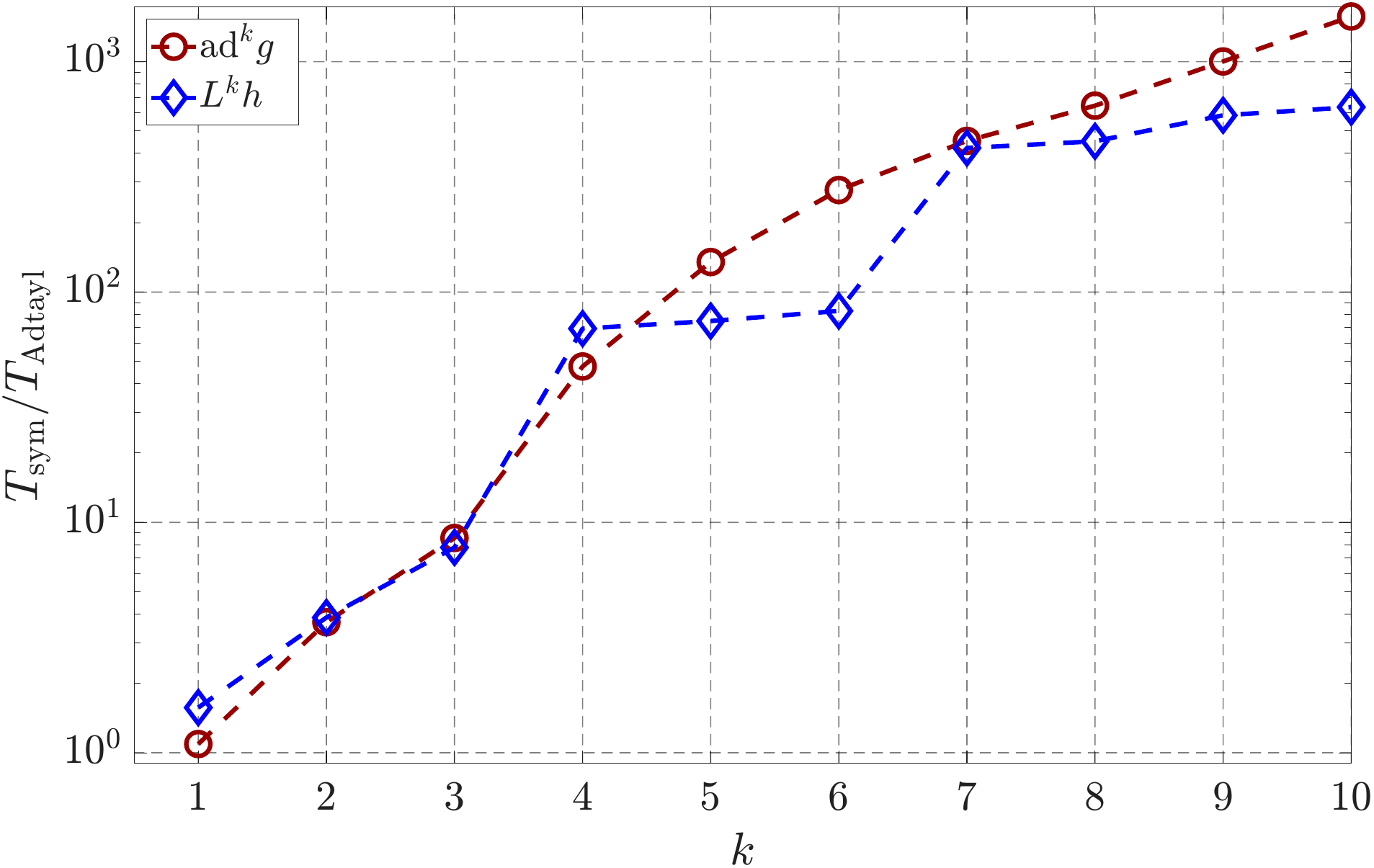}
    \caption{ Ratio of symbolic and \Adtayl  times.\label{fig:ratio}}
  \end{subfigure}
  \caption{Performance comparison for evaluation of $\adfk{k} g(x_0)$ and
  $\lider{k} h(x_0)$ using \matlab Symbolic Math Toolbox and \Adtayl.
  The time is shown on a logarithmic scale.}
  \label{fig:perf}
  \end{figure}
\begin{figure}[htb]
  \centering
  \includegraphics[width=\scl\linewidth]{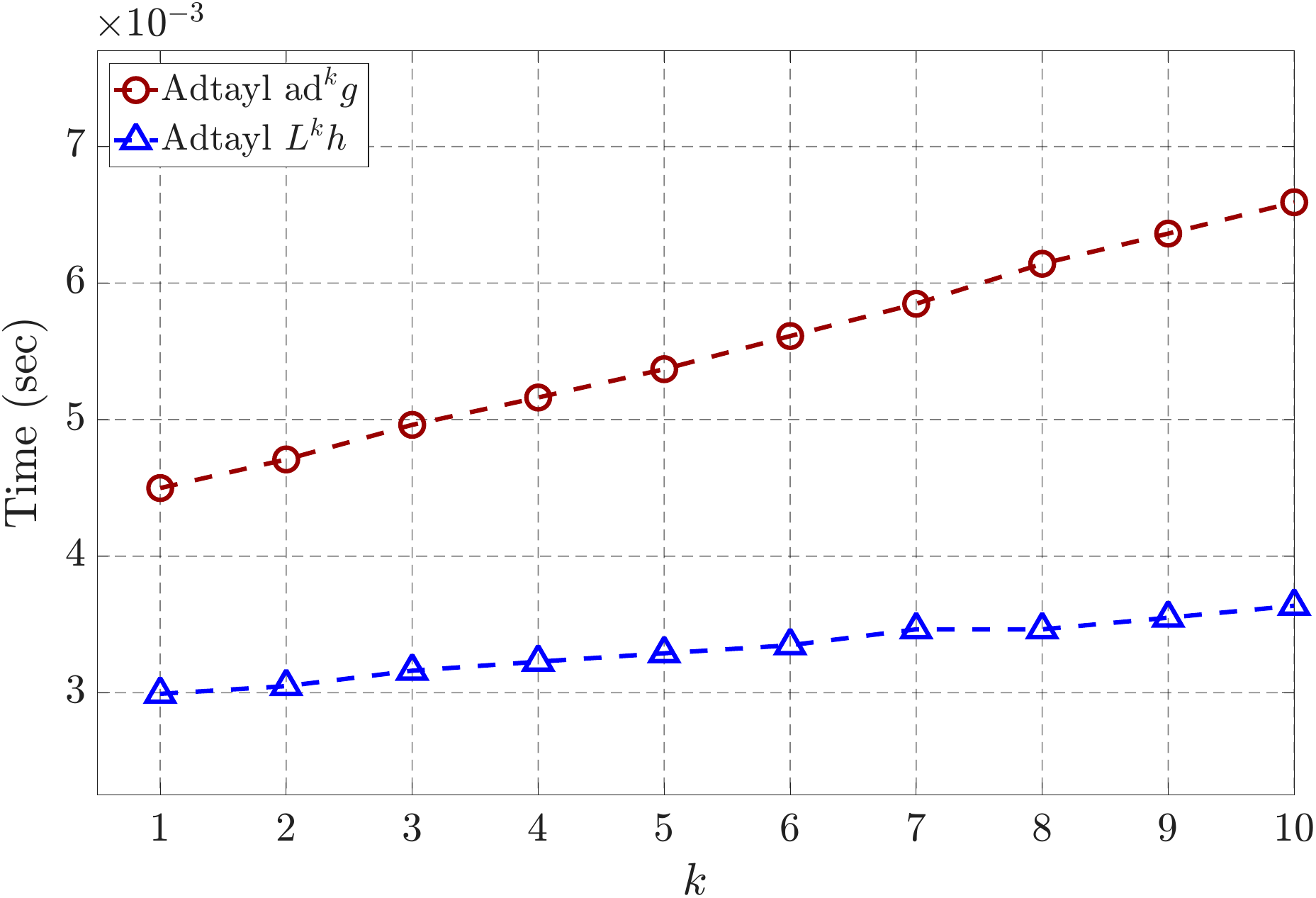}
  \caption{\Adtayl  time for computing $\adfk{k} g(x_0)$ and
  $\lider{k} h(x_0)$.}
  \label{fig:gantryAdtayl}
\end{figure}

\Cref{fig:combined} shows the time for evaluating 
$\adfk{k} g(x_0)$ and $\lider{k} h(x_0)$ versus order
$k$.  \Cref{fig:ratio} shows the ratio of symbolic and \Adtayl times\,
$T_{\mathrm{sym}}(k) / T_{\mathrm{Adtayl}}(k)$ at order $k$
for both $\adfk{k} g(x_0)$ and $\lider{k} h(x_0)$.
\Cref{fig:gantryAdtayl} shows the \Adtayl timings more clearly.  

The \Adtayl timings grow slowly with $k$ (approximately linearly over this
range), whereas the symbolic approach exhibits roughly exponential growth as
$k$ increases, due to the increasing size and complexity of intermediate
symbolic expressions.
For $k=10$, the symbolic cost is roughly 1000 times larger
than that of the \Adtayl method.

\begin{figure}[htb]
  \centering
 \includegraphics[width=\scl\linewidth]{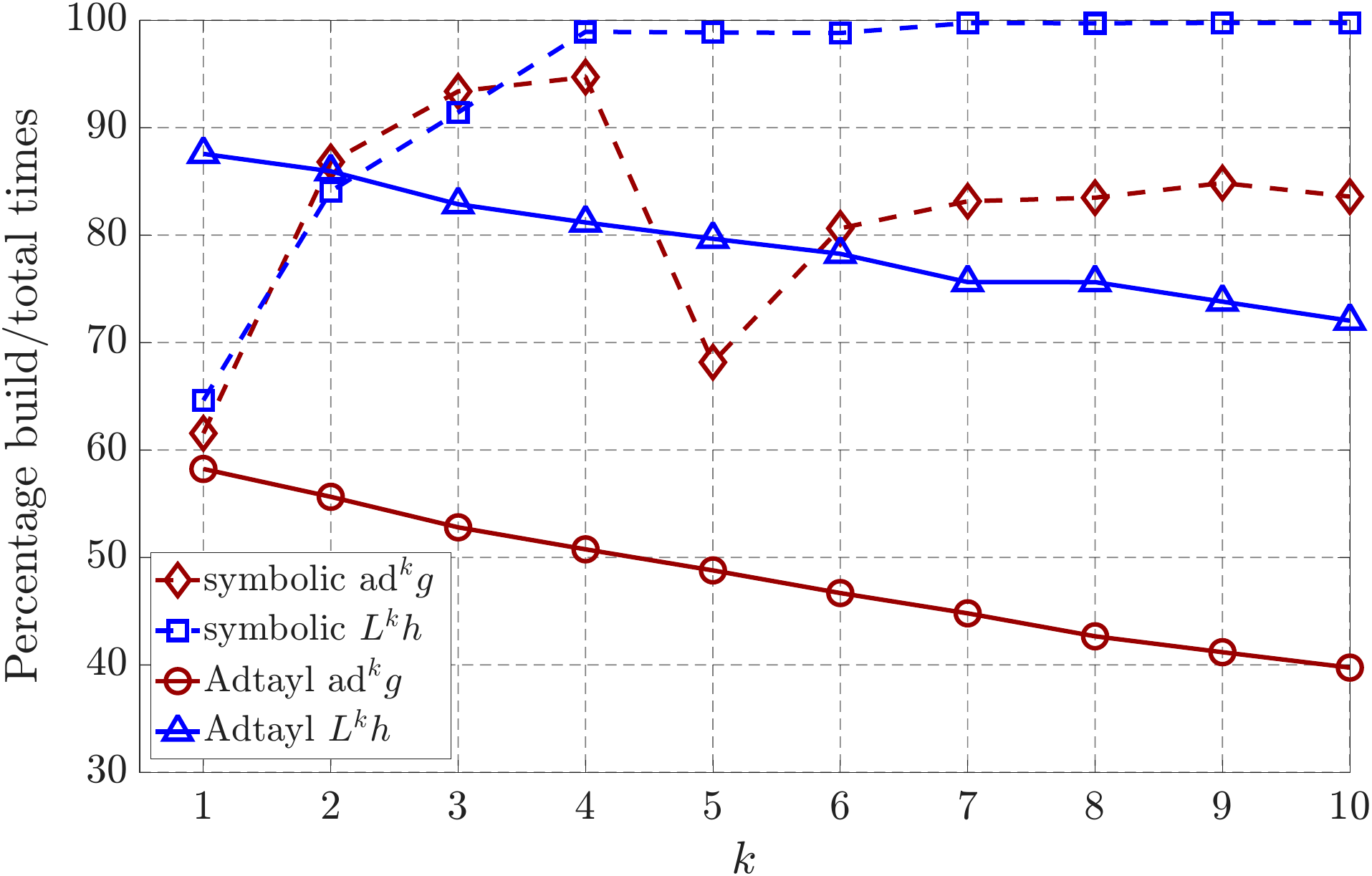}
\caption{Percentage of total runtime spent on symbolic expression
construction in the Symbolic Toolbox and on \cl generation in \Adtayl.\label{fig:build-fraction}}
\end{figure}

Figure \ref{fig:build-fraction} illustrates the percentage of total runtime dedicated to symbolic expression construction and \cl generation.
For the symbolic computation, the time for expression construction becomes increasingly dominant as the derivative order $k$ increases, approaching $100\%$ for $L^k h$ and fluctuating between $70\%$ and $90\%$ for $\adfk{k} g$ for $k \ge 4$.
In contrast, for \Adtayl, the time spent building the code list is fixed.
Consequently, its relative contribution to the total runtime decreases linearly as $k$ increases, cf.~Figure \ref{fig:gantryAdtayl}.

\section{Conclusions}\label{sc:concl}

\subsection{Software structure}

This article is primarily about computing \Lders of various kinds, but it also gives insights on numerical software package design.
Computing \Lders was an application suggested to us after \Adtayl was essentially written, and the economy of code shown in \cref{sc:lder} surprised us.

This economy is partly due to making \odets accept an ODE function $f$ with the same syntax used by the \matlab solvers, as far we could.
It is also due to our work to make \adtayl objects behave like native \matlab arrays as far as possible in terms of indexing, concatenating, slicing and basic matrix algebra.
Hence $f$, or its \cl, can evaluate on a numeric or on an \adtayl array; if the latter, either in one pass purely algebraically, or using the multi-pass recurrences needed for solving a differential equation.

In this way \Adtayl aspires towards the \term{differentiable programming} paradigm, which roughly speaking aims to make derivatives and differentiation---ordinary or partial, and with respect to variables, parameters or initial conditions---into first-class objects and operators of a computing platform.

\Adtayl does so imperfectly because individual \cl instructions are scalar.
Even the simple ODE $\xp=2x$ with $n$-vector $x$---i.e.\ $n$ independent equations $\xp_i=2x_i$---needs the $2x$ to become a loop in code for the $f$ of \rf{ivp}.
To ``vectorize'' the \cl would bring \Adtayl much further towards differentiable programming, but represents a major effort.

\subsection{Performance and numerics}

\subsubsection{Runtime scaling with order in \Adtayl.}
If the cost of computing derivatives up to order $k$ is dominated by arithmetic
operations (including array indexing), the total runtime is expected to scale
quadratically with $k$. However, for the order range considered here
($k \le 10$), the observed runtime increases approximately linearly with $k$.

The \cl is executed once per order.
For the range of $k$ considered here, profiling shows that the cost is dominated by order-independent overheads.
These include array accesses and scalar assignments used to decode each \cl entry (e.g.\ loading the opcode and operand indices), branching statements used to dispatch opcode executions, and special-case checks.
Although the algorithm also contains the convolution-type operations \rf{tcm}, \rf{tcd}, and \rf{ddot}, whose costs are proportional to $k$, their loops remain short for small $k$ and contribute  weakly to the total time.
For $k$ up to $10$, the cost of these convolutions grows linearly; quadratic behaviour (still mild) becomes apparent only at much higher orders---in double precision, for $k\gtrsim 50$.

\subsubsection{Numerical error.}
As for numerical error, we took as reference solution one computed by \matlab \li{vpa} in 40 digits of precision.
In double precision, we found \Adtayl's (normwise relative) error in $\lider{k} h$ grew very modestly as $k$ went from 1 to 10---from around roundoff to four times this.
The $\adfk{k} g$ error grew more, by a factor around 200.

For LIEDRIVERS, the errors for $\lider{k} h$ were equally good on this $k$ range.
However, we observed that its errors for $\adfk{k} g$ grew rapidly, by several orders of magnitude.
Since the difference is that the latter uses a Jacobian and the former doesn't, we wonder if there is a problem in the method used in ADOL-C to compute this Jacobian.

\subsection{Other work}

Following the theory described in \cite{robenack2005computation}, we can compute {\em observability matrices} compactly with \Adtayl, by code similar to that given above for \Lders.
We aim to report on this in future work.

\subsubsection*{Acknowledgments.}
The first author acknowledges the support of the Natural Sciences and Engineering
Research Council of Canada (NSERC), grant RGPIN-2019-07054.

\appendix 
\section{Gantry crane:  {\unboldmath $f$, $g$, and $h$} expressions\label{sc:gantryeqns}} 

With
\[
	x =
	\begin{pmatrix}
		x_1, x_2  , x_3, x_4
	\end{pmatrix}^\top
	=
	\begin{pmatrix}
		z ,
		\phi,
		\dot z ,
		\dot\phi
	\end{pmatrix}^\top,
\]
where $z$ is the position of the center of mass of the cart, $\phi$ is the 
angle of the cable measured from the vertical, and $(\dot z,\dot\phi)$ are the
corresponding velocities, the $f$ and $g$ in \rf{gantry} are

\[
	f(x) =
	\begin{pmatrix}
		x_3                  \\[0.5em]
		x_4                  \\[0.5em]
		\dfrac{m \ell x_4^{2} \sin x_2 {+} m G \sin x_2 \cos x_2}
		{m \sin^{2} x_2 {+}M} \\[1.5em]
		-\dfrac{m \ell x_4^{2} \sin x_2 \cos x_2 {+} (m {+}M)  G\sin x_2}
		{\ell\bigl(m \sin^{2} x_2 {+}M\bigr)}
	\end{pmatrix} \quad\text{and}\quad
	g(x) =
	\begin{pmatrix}
		0                             \\[0.5em]
		0                             \\[0.5em]
		\dfrac{1}{m \sin^{2} x_2 {+} M} \\[1.5em]
		\dfrac{{-}\cos x_2}{\ell\bigl(m \sin^{2} x_2 {+}M\bigr)}
	\end{pmatrix}
\]
The output map 
$
	h(x) =
	\begin{pmatrix}
			\ell\sin x_2 {+}x_1, 
			\ell\cos x_2
	\end{pmatrix}^\top
$ gives the Cartesian coordinates of the load.


\end{document}